\documentclass{article}
\usepackage[margin=1in]{geometry}
\usepackage{amsmath,amsthm,amssymb,mathtools, enumitem, setspace}
\usepackage[style=numeric, giveninits, url=false]{biblatex}
\usepackage[hidelinks]{hyperref}
\newtheorem{thm}{Theorem}[section]
\newtheorem{prop}[thm]{Proposition}
\newtheorem{cor}[thm]{Corollary}
\newtheorem{lem}[thm]{Lemma}
\theoremstyle{definition}
\newtheorem{defn}[thm]{Definition}
\newtheorem{example}[thm]{Example}
\newcommand{\N}{\mathbb{N}}
\DeclareMathOperator{\Ap}{Ap}
\title{Parity distributions among gaps of free numerical semigroups}
\author{
  Caleb M.\ Shor\thanks{\tt{caleb.shor@wne.edu}} \\ 
  Department of Mathematics \\ 
  Western New England University \\ 
  Springfield, MA, USA
}
\date{\today{}}

\begin{document}
\maketitle
\begin{abstract}
    In this paper, we extend recent results about the distribution of even and odd gaps of a numerical semigroup. We find that, for any numerical semigroup, the distribution can be computed in terms of the numbers of or the sums of odd and even elements in a corresponding Ap\'ery set. With free numerical semigroups specifically, we show that there are always at least as many odd gaps as even gaps, with equality precisely when the generating elements are all odd. We then specialize these results to the cases of numerical semigroups generated by compound and geometric sequences.
\end{abstract}

\section{Introduction}
In \cite{Cho03062025}, Cho, Lee, \& Nam investigated the distribution of odd and even gaps of numerical semigroups. They gave general results in the case where one knows the Ap\'ery set of a numerical semigroup relative to an odd generating element. For this, one needs only to know the difference between the numbers of even and odd elements of the Ap\'ery set. They also gave explicit results for a few families of numerical semigroups which have Ap\'ery sets that are known in the literature. Notably, in the case of a numerical semigroup $S$ generated by a compound sequence for which all minimal generating elements are odd, they showed that $S$ has as many odd gaps as it has even gaps. In their concluding remarks, they mentioned an open problem about computing the distribution of odd and even gaps in the case of a numerical semigroup generated by a compound sequence where at least one of the minimal generating elements is even.

In this paper, we expand on their results in two ways. First, with Proposition~\ref{prop:parity-difference} we extend their general results to the case where one knows the Ap\'ery set of a numerical semigroup relative to any nonzero (i.e., not just odd) element. Second, we consider the problem about the distribution of odd and even gaps for numerical semigroups generated by any compound sequence. We tackle this problem more generally, considering the same problem for so-called free numerical semigroups, of which those generated by compound sequences are a special case. We find a formula (Theorem~\ref{thm:main result}) for the difference between the numbers of odd and even gaps for any free numerical semigroup, concluding that there are always at least as many odd gaps as even gaps, with equality precisely when the minimal generating elements of the semigroup are all odd. 

\subsection{Organization}
This paper is organized as follows. In Section~\ref{sec:prelims}, we provide some background on numerical semigroups and their gaps along with the primary tools that we will use to understand them. This will culminate in Proposition~\ref{prop:parity-difference} which, for any numerical semigroup $S$, gives a formula for $O(G(S))-E(G(S))$, the difference between the number of odd gaps of $S$ and the number of even gaps of $S$. The formula is written in terms of an Ap\'ery set of $S$ relative to some nonzero element $t$ of $S$. We then illustrate this result with an example.

In Section~\ref{sec:free-numerical-semigroups}, we look specifically at the case of free numerical semigroups along with the sequences, known as telescopic sequences, that one uses to generate them. With knowledge of the Ap\'ery sets of free numerical semigroups, we use an identity from Gassert \& Shor \cite{GassertShor17} to compute the difference between the numbers of odd gaps and even gaps of any free numerical semigroup. Our main result is Theorem~\ref{thm:main result}. Immediately afterward, we have two corollaries: with one we give a simple formula for the difference in the case where all but one of the generating elements of a free numerical semigroup are even; and with the other we show that a free numerical semigroup $S$ has the same number of odd gaps as even gaps precisely when the minimal generating elements of $S$ are all odd. We then apply these results to numerical semigroups generated by compound sequences and geometric sequences, and we round off this section with some explicit examples.

Finally, in Section~\ref{sec:final-thoughts}, we describe two additional ways to determine the distribution of even and odd gaps of a numerical semigroup generated by a compound sequence when there are even generating elements. We close with a question for further study.

\section{Numerical semigroups}\label{sec:prelims}
\subsection{Preliminaries} 
For a comprehensive introduction to numerical semigroups, we refer the reader to \cite{RosalesGarciaSanchez09}.

Let $\N$ be the set of positive integers. Let $\N_0$ be the set of nonnegative integers, a monoid under addition. A \emph{numerical semigroup} $S$ is a submonoid of $\N_0$ with finite complement. Elements of this complement, denoted $G(S)=\N_0\setminus S$, are the \emph{gaps} of $S$.

For a sequence of positive integers $(n_1,\dots,n_k)$, the set of all nonnegative linear combinations of $n_1,\dots,n_k$ is denoted 
\[
    \langle n_1,\dots,n_k \rangle
    = \left\{\sum\limits_{i=1}^k c_in_i : c_1,c_2,\dots,c_k\in\N_0\right\}.
\]
It is known that $\langle n_1,\dots,n_k\rangle$ is a numerical semigroup if and only if $\gcd(n_1,\dots,n_k)=1$. Furthermore, every numerical semigroup $S$ can be written in the form $S=\langle n_1,\dots,n_k\rangle$ for some sequence of positive integers $(n_1,\dots,n_k)$ with $\gcd(n_1,\dots,n_k)=1$. We say that $S$ is minimally generated by $(n_1,\dots,n_k)$ if no proper subsequence generates $S$. It happens that every numerical semigroup $S$ has a unique minimal generating sequence. Such a minimal generating sequence is necessarily finite.

An important and useful tool for working with numerical semigroups is an Ap\'ery set.
\begin{defn}[Ap\'ery set]
    For $S$ a numerical semigroup and $t$ some nonzero element of $S$, the Ap\'ery set of $S$ relative to $t$ is
    \[\Ap(S;t)=\{s\in S : s-t\not\in S\}.\]
\end{defn}
Put another way, $\Ap(S;t)$ consists of exactly $t$ elements, each the smallest element of $S$ in its congruence class modulo $t$. Many properties of a numerical semigroup can be derived directly from one of its Ap\'ery sets. We will use the following identity from Gassert \& Shor \cite{GassertShor17}.

\begin{prop}[{\cite[Theorem 2.3]{GassertShor17}}]
    \label{prop:shor-identity}
    Let $S$ be a numerical semigroup with gap set $G(S)$. For any nonzero $t\in S$, and for any arithmetic function $f$,
    \[\sum\limits_{g\in G(S)}\left(f(g+t)-f(g)\right) 
    = \sum\limits_{a\in\Ap(S;\, t)}f(a)-\sum\limits_{k=0}^{t-1}f(k).\]
\end{prop}

As an example, using the function $f(n)=n$ with  Proposition~\ref{prop:shor-identity} results in a formula for $g(S)$, the \emph{genus} of $S$, which is the number of gaps of $S$:
\[g(S)=\frac{1-t}{2}+\frac{\sigma(\Ap(S;t))}{t},\] 
where $\sigma(T)$ denotes the sum of the elements of a finite set $T$. (This result is originally due to Selmer. See \cite[Section 2]{Selmer1977}.)

\subsection{Distribution of odd and even gaps given any Ap\'ery set}
In this paper, we are interested in the difference between the numbers of odd gaps of $S$ and even gaps of $S$, where $S$ is some numerical semigroup. For any finite set $T$ of integers, let $O(T)$ (resp., $E(T)$) denote the number of odd (resp., even) integers in $T$. Our goal is to understand the quantity $O(G(S))-E(G(S))$ without explicitly computing $G(S)$.

Cho, Lee, \& Nam give a formula (\cite[Theorem 3.1]{Cho03062025}) for $O(G(S))-E(G(S))$ in terms of $\Ap(S;t)$ in the case where $t$ is an odd minimal generating element of $S$. With Proposition~\ref{prop:shor-identity}, we are able to extend the result for any nonzero $t\in S$, odd or even. We also note that $t$ need not be a generating element of $S$ here. 

In what follows, for any finite set $T$ of integers, let $\sigma_O(T)$ (resp., $\sigma_E(T)$) denote the sum of the odd (resp., even) elements of $T$.

\begin{prop}\label{prop:parity-difference}
    For $S$ a numerical semigroup and some nonzero $t\in S$, let $\mathcal{A}_t:=\Ap(S;t)$. Then 
    \[O(G(S))-E(G(S)) = -\frac{1}{2} + 
    \begin{dcases*}
            \frac{E(\mathcal{A}_t)-O(\mathcal{A}_t)}{2} & if $t$ is odd, \\
            \frac{\sigma_O(\mathcal{A}_t)-\sigma_E(\mathcal{A}_t)}{t} & if $t$ is even.
    \end{dcases*}\]
\end{prop}
\begin{proof}
    If $t$ is odd, we use $f(n)=(-1)^n$ with Proposition~\ref{prop:shor-identity}:
    \[
        \sum\limits_{g\in G(S)}\left((-1)^{g+t}-(-1)^g\right) 
        = \sum\limits_{a\in \mathcal{A}_t}(-1)^a - \sum\limits_{k=0}^{t-1}(-1)^k.
    \]
    Since $t$ is odd, $(-1)^{g+t}=-(-1)^g$ and $\sum_{k=0}^{t-1}(-1)^t=1$. Thus,
    \[
        O(G(S))-E(G(S)) 
        = -\sum\limits_{g\in G(S)}(-1)^g 
        = \frac{1}{2}\left(-1+\sum\limits_{a\in A_t}(-1)^a\right),\]
    from which the given formula follows.

    If $t$ is even, the function $f(n)=(-1)^n$ doesn't help. We instead use $f(n)=n(-1)^n$ with Proposition~\ref{prop:shor-identity}: 
    \[
        \sum\limits_{g\in G(S)}\left((g+t)(-1)^{g+t}-g(-1)^g\right) 
        = \sum\limits_{a\in \mathcal{A}_t}a(-1)^a 
        - \sum\limits_{k=0}^{t-1}k(-1)^k.
    \]
    Since $t$ is even, $(-1)^{g+t}=(-1)^g$ and $\sum_{k=0}^{t-1}k(-1)^k=-t/2$. Thus,
    \[O(G(S))-E(G(S)) 
    = -\sum\limits_{g\in G(S)}(-1)^g 
    = -\frac{1}{t}\left(\frac{t}{2} + \sum\limits_{a\in \mathcal{A}_t}a(-1)^a \right),\]
    from which the given formula follows.
\end{proof}

As a result, in order to compute $O(G(S))-E(G(S))$ for any numerical semigroup $S$, we just need a nice description of $\Ap(S;t)$ for some nonzero $t\in S$. We will illustrate this with a small example.
\begin{example}\label{ex:first-ex}
    Suppose $S=\langle 6,15,20\rangle$. We will compute $O(G(S))-E(G(S))$ using Proposition~\ref{prop:parity-difference} with $\Ap(S;6)$ and with $\Ap(S;15)$. We will then explicitly determine the gap set $G(S)$ to confirm our results.

    First, let $\mathcal{A}_{15}=\Ap(S;15)$. With some computation, we find that 
    \[
        \mathcal{A}_{15}=\Ap(S;15)=\{0,6,12,18,20,24,26,32,38,40,44,46,52,58,64\},
    \]
    a set consisting of $O(\mathcal{A}_{15})=0$ odd integers and $E(\mathcal{A}_{15})=15$ even integers. By Proposition~\ref{prop:parity-difference}, since 15 is odd, we get
    \[
        O(G(S))-E(G(S)) 
        = -\frac{1}{2}+\frac{E(\mathcal{A}_{15})-O(\mathcal{A}_{15})}{2} 
        = -\frac{1}{2}+\frac{15-0}{2}
        =7.
    \]
    
    Next, let $\mathcal{A}_6=\Ap(S;6)$. With some computation, we find that $\mathcal{A}_6=\Ap(S,6)=\{0,15,20,35,40,55\}$. The sum of the odd elements of $\mathcal{A}_6$ is $\sigma_O(\mathcal{A}_6)=15+35+55=105$. The sum of the even elements of $\mathcal{A}_6$ is $\sigma_E(\mathcal{A}_6)=0+20+40=60$. By Proposition~\ref{prop:parity-difference}, since 6 is even, we get
    \[
        O(G(S))-E(G(S)) 
        = -\frac{1}{2}+\frac{\sigma_O(\mathcal{A}_6)-\sigma_E(\mathcal{A}_6)}{6}
        =-\frac{1}{2}+\frac{105-60}{6}
        = 7.
    \]

    To verify our results, some more computation reveals the gap set of $S$ to be 
    \[
        G(S)=\{1,2,3,4,5,7,8,9,10,11,13,14,16,17,19,22,23,25,28,29,31,34,37,43,49\},
    \]
    consisting of 16 odd gaps and 9 even gaps. We can verify directly that $O(G(S))-E(G(S))=16-9=7$.
\end{example}

\section{Specialization to free numerical semigroups}\label{sec:free-numerical-semigroups}
\subsection{Telescopic sequences and free numerical semigroups}
We are now ready to specialize our results to free numerical semigroups. For a more detailed treatment, see \cite{Shor19a} and the references therein. 

We begin by defining telescopic sequences (also called smooth sequences).
\begin{defn}[Telescopic sequence]
    For some $k\in\N_0$, let $T=(t_0,\dots,t_k)$ be a sequence of positive integers. For $0\le i\le k$, let $T_i=(t_0,\dots,t_i)$, and let $d_i=\gcd(T_i)$. Then, let $c_j=d_{j-1}/d_j$ for all $1\le j\le k$, and define the sequence $c(T)=(c_1,\dots,c_k)$. If $c_jt_j\in\langle T_{j-1}\rangle$ for all $j=1,\dots,k$, then $T$ is a \emph{telescopic} sequence.
\end{defn}
There are many families of sequences that are telescopic. For instance, any two-term $(a,b)$ sequence is telescopic, as is any geometric sequence $(a^k, a^{k-1}b, \dots,b^k)$. As we will see in  Section~\ref{sec:compound}, compound sequences are telescopic as well.

\begin{lem}\label{lem:t0productofci}
    Suppose $T=(t_0,\dots,t_k)$ is a telescopic sequence with $c(T)=(c_1,\dots,c_k)$. If $\gcd(T)=1$, then for any $1\le j\le k$, we have $\prod\limits_{i=j}^k c_i=d_{j-1}$. In particular, we have $\prod\limits_{i=1}^k c_i=t_0$.
\end{lem}
\begin{proof}
    Suppose $1\le j\le k$. Then
    \[
        \prod\limits_{i=j}^k c_i
        = c_j c_{j+1}\cdots c_k
        = \dfrac{d_{j-1}}{d_j}\dfrac{d_j}{d_{j+1}}\cdots \dfrac{d_{k-1}}{d_k} = \dfrac{d_{j-1}}{d_k}.
    \]
    Since $d_k=\gcd(T_k)=\gcd(T)=1$, we conclude that $\prod\limits_{i=j}^k c_i = d_{j-1}$.

    In particular, with $j=1$, we get
    \[
        \prod\limits_{i=1}^k c_i = d_0 = \gcd(T_0)=\gcd(t_0)=t_0.\qedhere
    \]
\end{proof}
\begin{defn}[Free numerical semigroup]
    Let $T$ be a telescopic sequence. If $\gcd(T)=1$, then  $S=\langle T\rangle$ is called a \emph{free numerical semigroup}.
\end{defn}

Free numerical semigroups are nice to work with because they have easily describable Ap\'ery sets.

\begin{prop}[{\cite[Proposition 6]{GassertShor17}}]
    \label{prop:free-apery-description}
    Suppose $S$ is a free numerical semigroup. Then $S=\langle T\rangle$ for some telescopic sequence $T=(t_0,\dots,t_k)$ with $\gcd(T)=1$. Let $c(T)=(c_1,\dots,c_k)$. Then
    \[
        \Ap(S;t_0)
        =\left\{\sum\limits_{i=1}^k n_it_i : 0\le n_i<c_i\text{ for } i=1,\dots,k\right\}.
    \]
\end{prop}

\subsection{Distribution of odd and even gaps of a free numerical semigroup}
Given our description of an Ap\'ery set of a free numerical semigroup $S$, we can now use Proposition~\ref{prop:parity-difference} to investigate the parity distributions among gaps of $S$. This is our main result.
\begin{thm}\label{thm:main result}
    Suppose $S$ is a free numerical semigroup generated by a telescopic sequence $T=(t_0,\dots,t_k)$. Let $c(T)=(c_1,\dots,c_k)$, and let $c_0=1$. Let $m=\min\{0\le i\le k : t_i\text{ is odd}\}$ and $I=\{0\le i\le k : t_i\text{ is even}\}$. Then
    \[
        O(G(S))-E(G(S)) 
        = -\frac{1}{2}+\frac{c_mt_m}{2t_0}\prod\limits_{i\in I}c_i.
    \]
    In particular, if $t_0$ is odd, then 
    \[
        O(G(S))-E(G(S)) = -\frac{1}{2}+\frac{1}{2}\prod\limits_{i\in I}c_i.
    \]
\end{thm}
\begin{proof}
    To start, note that $m$ is well-defined, for otherwise $\gcd(T)\ne1$ and $S$ is not a numerical semigroup. In addition to the notation defined in the statement of the theorem, we will let $\mathcal{A}_{t}=\Ap(S;t)$. In what follows, we will consider the cases of $t_0$ odd and even separately.

    First, suppose $t_0$ is odd. In this case, $m=0$. Since $c_0=1$ and $t_m=t_0$, our goal is to show that $O(G(S))-E(G(S))=-1/2 + 1/2\prod\limits_{i\in I}c_i$.
    
    By Proposition~\ref{prop:parity-difference}, we need to compute $E(\mathcal{A}_{t_0})-O(\mathcal{A}_{t_0})$. With our description of $\mathcal{A}_{t_0}$ from Proposition~\ref{prop:free-apery-description}, we have
    \begin{align*}
        E(\mathcal{A}_{t_0})-O(\mathcal{A}_{t_0})
        &= \sum\limits_{a\in \mathcal{A}_{t_0}} (-1)^a \\
        &= \sum\limits_{n_1=0}^{c_1-1} \cdots \sum\limits_{n_k=0}^{c_k-1} (-1)^{n_1t_1+\cdots+n_kt_k} \\
        &= \prod\limits_{i=1}^k \sum\limits_{n_i=0}^{c_i-1}(-1)^{n_it_i}.
    \end{align*}
    Observe that for any positive integers $c$ and $t$,
    \begin{equation}\label{eqn:summation-cases}
        \sum\limits_{n=0}^{c-1}(-1)^{n t} = 
        \begin{dcases*}
            c & if $t$ is even, \\
            1 & if $t$ is odd and $c$ is odd, \\
            0 & if $t$ is odd and $c$ is even.
        \end{dcases*}
    \end{equation}
    By Lemma~\ref{lem:t0productofci}, $t_0=c_1\cdots c_k$. Since $t_0$ is odd, all $c_i$ terms are odd. The summation in Eq.~\eqref{eqn:summation-cases} will never be 0. In particular, if we let $I_0=\{1\le i\le k : t_i\text{ is even}\}$, then we have 
    \[
        E(\mathcal{A}_{t_0})-O(\mathcal{A}_{t_0}) 
        = \prod\limits_{i\in I_0}c_i.
    \]
    We defined $I$ in the statement of the theorem to consider indices $0\le i\le k$. Since $t_0$ is odd, $0\not\in I$ and hence $I_0=I$. By Proposition~\ref{prop:parity-difference}, we conclude that
    \[
        O(G(S))-E(G(S)) 
        = -\frac{1}{2}+\frac{E(\mathcal{A}_{t_0})-O(\mathcal{A}_{t_0})}{2} 
        = -\frac{1}{2}+\frac{1}{2}\prod\limits_{i\in I}c_i,
    \]
    as desired.

    Next, suppose $t_0$ is even. With our definition of $m$ as the index of the first odd $t_i$, we can draw a few conclusions. First, $t_0,\dots,t_{m-1}$ are even, and hence $d_0,\dots,d_{m-1}$ are even as well. Next, $d_m$ is odd, and this means that $c_m$ is even and $c_{m+1},\dots,c_k$ are odd. In particular, $i=m$ is the unique index among $0\le i\le k$ for which $t_i$ is odd and $c_i$ is even.

    Now, in order to use Proposition~\ref{prop:parity-difference}, we compute 
    \begin{align*}
        \sigma_E(\mathcal{A}_{t_0})-\sigma_O(\mathcal{A}_{t_0})
        &= \sum\limits_{a\in \mathcal{A}_{t_0}}a(-1)^a \\
        &= \sum\limits_{n_1=0}^{c_1-1}\cdots \sum\limits_{n_k=0}^{c_k-1} (n_1t_1+\cdots+n_kt_k)(-1)^{n_1t_1+\cdots+n_kt_k} \\
        &= \sum\limits_{j=1}^k F_j,
    \end{align*}
    where, for $j=1,\dots,k$, we have
    \[
        F_j 
        = \left( 
            \prod\limits_{\substack{i=1, \\ i\ne j}}^k \sum\limits_{n_i=0}^{c_i-1}(-1)^{n_it_i} 
        \right) \cdot
        \sum\limits_{n_j=0}^{c_j-1}n_jt_j(-1)^{n_jt_j}.
    \]
    We saw the summation over $n_i$ above earlier in Eq.~\eqref{eqn:summation-cases}. Since the case of $t_i$ odd and $c_i$ even occurs exactly once, for $i=m$, we conclude that $F_j=0$ for all $j\ne m$.

    Let $I_m = \{1\le i\le k, i\ne m : t_i\text{ is even}\}$. To compute $F_m$, we note that 
    \[
        \sum\limits_{n_m=0}^{c_m-1}n_mt_m(-1)^{n_mt_m} = 0 - t_m + 2t_m - \cdots - (c_m-1)t_m = -\frac{c_mt_m}{2}.
    \]
    Thus,
    \[
        \sigma_E(\mathcal{A}_{t_0})-\sigma_O(\mathcal{A}_{t_0}) 
        = F_m 
        = -\frac{c_mt_m}{2}\prod\limits_{i\in I_m}c_i.
    \]
    Next, observe that since $t_m$ is odd, we have $m\not\in I$ (for $I$ as defined in the statement of the theorem), and hence $I_m=I$. Finally, by Proposition~\ref{prop:parity-difference},
    \[
        O(G(S))-E(G(S)) 
        = -\frac{1}{2}-\frac{\sigma_E(\mathcal{A}_{t_0})-\sigma_O(\mathcal{A}_{t_0})}{t_0} 
        = -\frac{1}{2}+\frac{c_mt_m}{2t_0}\prod\limits_{i\in I}c_i,
    \]
    as desired.
\end{proof}

We obtain a relatively simple result in the case where all but one of the generating elements in a telescopic sequence are even.
\begin{cor}\label{cor:all-but-one-generating-elts-even}
    Suppose $S$ is a free numerical semigroup generated by a telescopic sequence $T=(t_0,\dots,t_k)$. Further assume that there is exactly one odd generating element $t_m$. Then 
    \[O(G(S))-E(G(S)) = \frac{t_m-1}{2}.\]
\end{cor}
\begin{proof}
    If $t_m$ is odd and all other $t_i$ terms are even, then $I=\{0,1,\dots,k\}\setminus\{m\}$. Since $t_0=c_0c_1\cdots c_k$, we have $\prod\limits_{i\in I} c_i = t_0/c_m$. By Theorem~\ref{thm:main result}, 
    \[
        O(G(S))-E(G(S)) 
        = -\frac{1}{2}+\frac{c_mt_m}{2t_0}\cdot \frac{t_0}{c_m},
    \]
    and the result follows.
\end{proof}

We now show that there are always at least as many odd gaps as even gaps for any free numerical semigroup.

\begin{cor}\label{cor:more-odds-than-evens}
    Let $S$ be a free numerical semigroup generated by a minimal telescopic sequence $T=(t_0,\dots,t_k)$. Then there are at least as many odd gaps of $S$ as there are even gaps of $S$, with equality if and only if all terms in $T$ are odd.
\end{cor}
\begin{proof}
    Consider the sequence $c(T)=(c_1,\dots,c_k)$. Since $T$ is telescopic, $c_it_i\in\langle t_0,\dots,t_{i-1}\rangle$ for all $i=1,\dots,k$. And since $T$ is minimal, we may conclude that $t_i\not\in\langle t_0,\dots,t_{i-1}\rangle$ for all $i=1,\dots,k$. Hence, we have that $c_i>1$ for all $i=1,\dots,k$.
    
    Recall that $m=\min\{0\le i\le k : t_i\text{ is odd}\}$ and $I=\{0\le i\le k : t_i\text{ is even}\}$. We will consider $m=0$ and $m>0$ separately.
    
    By Theorem~\ref{thm:main result}, if $m=0$ then 
    \[
        O(G(S))-E(G(S))
        = -\frac{1}{2}+\frac{1}{2}\prod\limits_{i\in I}c_i.
    \]
    Hence, $O(G(S))\ge E(G(S))$, with equality if and only if $\prod_{i\in I}c_i=1$. Since the $c_i$ terms are all greater than 1, $O(G(S))=E(G(S))$ if and only if $I=\emptyset$, which means all terms in $T$ are odd.
    
    Now, suppose $m>0$, which implies $\{0,\dots,m-1\}\subseteq I$. Since we have a positive number of even elements in $T$, we need to show the number of odd gaps of $S$ is strictly larger than the number of even gaps of $S$. Appealing to Theorem~\ref{thm:main result} once again, we have 
    \[
        O(G(S)) - E(G(S))
        = -\frac{1}{2}+\frac{c_mt_m}{2t_0}\prod\limits_{i\in I}c_i.
    \]
    In order to show $O(G(S))>E(G(S))$, it will suffice to show that 
    \[
        c_mt_m\prod\limits_{i\in I}c_i
        > t_0.
    \]
    
    To start, we have by Lemma~\ref{lem:t0productofci} that 
    \[
        \prod\limits_{i=m+1}^k c_i 
        = d_m 
        = \gcd(t_0,\dots,t_m).
    \]
    Thus, $c_{m+1}\cdots c_k\le t_m$. If we have equality, then $t_m=d_m$. Since $m>0$, this implies $t_m\mid t_0,t_1,\dots,t_{m-1}$, contradicting the minimality of $T$. Hence, we must have $c_{m+1}\cdots c_k<t_m$.

    Multiplying both sides of this inequality by $c_1\cdots c_m$ (which is necessarily positive), we get $c_1\cdots c_k<t_mc_1\cdots c_m$. By Lemma~\ref{lem:t0productofci}, we have $t_0<t_m c_1\cdots c_m$. Since $\{0,\dots,m-1\}\subseteq I$, we conclude that 
    \[
        c_mt_m\prod\limits_{i\in I}c_i 
        \ge c_mt_m c_1\cdots c_{m-1} 
        > t_0,
    \]
    as desired.
\end{proof}
\subsection{Specialization to compound sequences}\label{sec:compound}
We now focus on numerical semigroups generated by compound sequences, which were first described in \cite{KiersONeillPonomarenko16}.
\begin{defn}[Compound sequence]
    For some $p\in\N_0$, let $A,B\in\N^p$ be sequences of positive integers. For $A=(a_1,\dots,a_p)$ and $B=(b_1,\dots,b_p)$, if $\gcd(a_i,b_j)=1$ for all $i\ge j$, then we say $(A,B)$ is a \emph{suitable pair of sequences}. Given a suitable pair $(A,B)$, we form the \emph{compound sequence} $C(A,B):=(n_0,\dots,n_p)\in\mathbb{N}^{p+1}$ defined by $n_0=a_1a_2\cdots a_p$ and, for $i\ge1$, $n_i=n_{i-1} b_i / a_i$. 
\end{defn} 
Note that the definition of a compound sequence in \cite{KiersONeillPonomarenko16} specifies that $2\le a_i<b_i$ for all $i$. With this condition in place, one has that $(n_0,\dots,n_p)$ is an increasing sequence. This is unnecessary for our work so we make no such restriction. The condition that $\gcd(a_i,b_j)=1$ for all $i\ge j$ is sufficient for our purposes.

The following proposition, from \cite{GassertShor16}, says that the condition where $a_i,b_i\ge2$ for all $i$ gives rise to a minimal generating sequence.

\begin{prop}{\cite[Propositions 2.6 and 2.7]{GassertShor16}}
    \label{prop:compound-gcd-minimal}
    Suppose $(A,B)$ is a suitable pair with $C(A,B)=(n_0,\dots,n_p)$. Then $\gcd(n_0,\dots,n_p)=1$, implying that $\langle n_0,\dots,n_p\rangle$ is a numerical semigroup. If we further have that $a_i,b_i\ge2$ for all $i$, then $(n_0,\dots,n_p)$ is a minimal generating sequence for this numerical semigroup.
\end{prop}

A compound sequence is a generalization of a geometric sequence, which one obtains in the case where $A=(a,a,\dots,a)$ and $B=(b,b,\dots,b)$. It turns out that compound sequences are telescopic.

\begin{prop}[{\cite[Corollary 3.11]{GassertShor17}}]\label{prop:compound-is-telescopic}
    If $CS=(n_0,\dots,n_p)$ is a compound sequence, then $CS$ is a telescopic sequence. In particular, if $CS=C(A,B)$ for a suitable pair $(A,B)$ with $A=(a_1,\dots,a_p)$ and $B=(b_1,\dots,b_p)$, then $c(CS)=A=(a_1,\dots,a_p)$.
\end{prop}

As a result, we can compute $O(G(S))-E(G(S))$ for a numerical semigroup generated by a compound sequence via Theorem~\ref{thm:main result}. For notation, we define the functions $m_2$ and $M_2$ of a sequence $Z=(z_1,\dots,z_p)$ of integers as follows:
\[
    m_2(Z):=\min\left(\{1\le i\le p : z_i \text{ is even}\}\cup\{p+1\}\right)
\]
and
\[
    M_2(Z):=\max\left(\{0\}\cup\{1\le i\le k : z_i\text{ is even}\} \right).
\]
Note that each function returns an element out of the index range if no such even integer exists in the sequence.

\begin{cor}\label{cor:distribution-for-compound}
    Suppose $A=(a_1,\dots,a_p)$ and $B=(b_1,\dots,b_p)$ are in $\N^p$, that $(A,B)$ is a suitable pair, and that $S=\langle C(A,B)\rangle$. Then
    \[
        O(G(S))-E(G(S)) = \dfrac{\alpha\beta-1}{2}
    \]
    where 
    \[
        \alpha = \prod\limits_{i=m_2(B)}^{p} a_i
        \quad\text{ and }\quad 
        \beta = \prod\limits_{i=1}^{M_2(A)} b_i.
    \]
\end{cor}
\begin{proof}
    By Proposition~\ref{prop:compound-is-telescopic}, $S$ is telescopic so we may use Theorem~\ref{thm:main result}. Our goal is to get a description of $I$, the set of indices of even generating elements of $S$.

    Since $C(A,B)=(n_0,\dots,n_p)$, we have that $n_m$ is the first odd generating element. If $m=0$, then $n_0=a_1\cdots a_p$ is odd, meaning all terms in $A$ are odd and thus $M_2(A)=0=m$. If $m>0$, then $n_m=b_1\cdots b_m a_{m+1}\cdots a_p$ is odd and $a_m$ is even. In other words, $M_2(A)=m$ in this case as well. Since $a_m$ is a factor of $n_0,\dots,n_{m-1}$, these generating elements are all even.

    Which other generating elements are even? Focus on $m_2(B)$. This gives the index of the first even term in $B$, and this term is a factor of $n_{m_2(B)}, \dots, n_p$. These generating elements are therefore even.

    Finally, note that since $(A,B)$ is a suitable pair, we have $M_2(A)<m_2(B)$, and we conclude that the generating elements $n_{M_2(A)},\dots,n_{m_2(B)-1}$ are odd. We can now describe $I$ explicitly:
    \[
        I
        =\{0,1,\dots,M_2(A)-1\}\cup\{m_2(B),\dots,p\}.
    \]

    By Theorem~\ref{thm:main result}, 
    \begin{align*}
        O(G(S))-E(G(S)) 
        &=-\frac{1}{2}+\frac{a_m n_m}{2n_0} \prod\limits_{i\in I}c_i \\
        &= -\frac{1}{2}+\dfrac{a_{M_2(A)}\cdot b_1\cdots b_{M_2(A)} a_{M_2(A)+1}\cdots a_p}{2 a_1\cdots a_p} \cdot a_1\cdots a_{M_2(A)-1}\cdot a_{m_2(B)}\cdots a_p.
    \end{align*}
    Canceling $a_1\cdots a_p$ from top and bottom, we obtain
    \[
        O(G(S))-E(G(S))
        = -\frac{1}{2}+\frac{b_1\cdots b_{M_2(A)} a_{m_2(B)}\cdots a_p}{2} 
        = -\frac{1}{2} + \frac{\alpha\beta}{2},
    \]
    as desired.
\end{proof}
In Corollary~\ref{cor:more-odds-than-evens}, we saw that any free numerical semigroup has at least as many odd gaps as even gaps, with equality precisely when the minimal generating sequence consists of odd terms. In the case of a numerical semigroup generated by a compound sequence, we obtain the following.
\begin{cor}
    Suppose $(A,B)$ is a suitable pair of sequences of integers that are all greater than 1 so that the compound sequence $C(A,B)$ generates the numerical semigroup $S$. Then $S$ has at least as many odd gaps as even gaps. Furthermore, the following are equivalent.
    \begin{enumerate}
        \item $O(G(S))=E(G(S))$.
        \item All terms in the compound sequence $C(A,B)$ are odd.
        \item All terms in the sequences $A$ and $B$ are odd.
    \end{enumerate}
\end{cor}

\subsection{Specialization to geometric sequences}
Consider the geometric sequence $R=(a^k, a^{k-1}b, \dots, b^k)$ for $a,b,k\in\N$. Then $\gcd(R)=1$ if and only if $\gcd(a,b)=1$. When this occurs, we say $S=\langle R\rangle$ is a \emph{geometric numerical semigroup}. As noted earlier, any geometric sequence is a special case of a compound sequence, which in turn is a special case of a free numerical semigroup. In the sequence $(a^k, a^{k-1}b, \dots, b^k)$, either all terms are odd, or all terms but one are even. By Corollaries~\ref{cor:all-but-one-generating-elts-even} and \ref{cor:more-odds-than-evens}, we immediately have the following result.

\begin{cor}
    Suppose $S$ is a geometric numerical semigroup, so $S=\langle a^k, a^{k-1}b,\dots, b^k\rangle$ for $a,b,k\in\N$ with $\gcd(a,b)=1$. Then at least one of $a$ and $b$ is odd, and we have 
    \[
        O(G(S))-E(G(S)) = 
        \begin{dcases*} 
            0 & if both $a$ and $b$ are odd, \\
            \frac{a^k-1}{2} & if $a$ is odd and $b$ is even, \\
            \frac{b^k-1}{2} & if $a$ is even and $b$ is odd.
        \end{dcases*}
    \]
\end{cor}
Just as we saw earlier, for $S$ a geometric numerical semigroup, we have $O(G(S))=E(G(S))$ if and only if the generating elements of $S$ are odd. Otherwise $S$ has more odd gaps than even gaps.

\subsection{Examples}
We conclude this section with a few examples to illustrate our results, beginning with a free numerical semigroup.

\begin{example}
    Consider the sequence $T=(t_0,t_1,t_2,t_3,t_4)=(120, 180, 100, 55, 22)$. To see if this sequence is telescopic, we compute $d_0=120$, $d_1=60$, $d_2=20$, $d_3=5$, $d_4=1$. Then $c_1=2$, $c_2=3$, $c_3=4$, and $c_4=5$. One can verify that $c_it_i\in\langle t_0,\dots,t_{i-1}\rangle$ for $i=1,2,3,4$. Hence, $T$ is telescopic. Furthermore, $\gcd(T)=1$, so we have a free numerical semigroup $S=\langle 120,180,100,55,22\rangle$. In the notation of Theorem~\ref{thm:main result}, we have $m=3$ and $I=\{0,1,2,4\}$. Then
    \[
        O(G(S))-E(G(S)) 
        = -\frac{1}{2}+\frac{c_3t_3}{2t_0}\prod\limits_{i\in I}c_i 
        = -\frac{1}{2}+\frac{4\cdot55}{2\cdot120}\cdot 1\cdot2\cdot3\cdot5
        =27.
    \]

    Note that this can also be handled by Corollary~\ref{cor:all-but-one-generating-elts-even} since $T$ is telescopic and all but one of the generating elements are even. Since the odd generating element is 55, the result is simply $(55-1)/2$.
\end{example}

We can now revisit Example~\ref{ex:first-ex}, in which we worked with a numerical semigroup generated by a compound sequence.

\begin{example}
    Let $A=(2,3)$ and $B=(5,4)$. Since $\gcd(a_i,b_j)=1$ for all $i\ge j$, $(A,B)$ is a suitable pair with corresponding compound sequence $C(A,B)=(2\cdot3,5\cdot3,5\cdot4)=(6,15,20)$ and numerical semigroup $S=\langle 6,15,20\rangle$. The even terms in $A$ and $B$ are $a_1$ and $b_2$. Hence, $M_2(A)=1$ and $m_2(B)=2$. We find 
    \[
        \alpha=\prod_{i=m_2(B)}^2a_i=a_2=3
        \quad\text{ and }\quad 
        \beta=\prod\limits_{i=1}^{M_2(A)}b_i=b_1=5.
    \]
    By Corollary~\ref{cor:distribution-for-compound}, we find that
    \[
        O(G(S))-E(G(S)) 
        = \dfrac{\alpha\beta-1}{2}
        =\dfrac{3\cdot5-1}{2}
        =7,
    \]
    agreeing with our work in Example~\ref{ex:first-ex}.
\end{example}

\begin{example}
    Let $A=(2,3,4,3,3,3)$ and $B=(5,5,5,5,8,5)$, a suitable pair. We obtain the compound sequence $C(A,B)$ and corresponding numerical semigroup 
    \[
        S
        =\langle C(A,B)\rangle
        =\langle 648, 1620, 2700, 3375, 5625, 15000, 25000\rangle.
    \]
    The maximal index of an even term in $A$ is $M_2(A)=3$ and the minimal index of an even term in $B$ is $m_2(B)=5$. Thus, 
    \[
        \alpha
        =\prod\limits_{i=m_2(B)}^6a_i
        =a_5a_6
        =3\cdot3
        =9
        \quad\text{ and }\quad 
        \beta
        =\prod\limits_{i=1}^{M_2(A)}b_i
        = b_1b_2b_3
        = 5\cdot5\cdot5,
    \]
    from which we conclude
    \[O(G(S))-E(G(S)) = \dfrac{\alpha\beta-1}{2}=\dfrac{a_5a_6b_1b_2b_3-1}{2}=562.\]
\end{example}

Finally, we note that we never assumed our generating sequences were minimal in Theorem~\ref{thm:main result} or Corollary~\ref{cor:distribution-for-compound}. Indeed, our results hold for non-minimal generating sequences. We demonstrate this with an example.

\begin{example}
    Let $A=(6,5)$ and $B=(1,2)$. Then $(A,B)$ is a suitable pair with $C(A,B)=(30,5,2)$. Of course, this is not a minimal sequence due to the 1 in $B$. We have $S=\langle C(A,B)\rangle=\langle 30,5,2\rangle=\langle 5,2\rangle$, but we may proceed anyway. We find $M_2(A)=1$ and $m_2(B)=2$. Thus, $\alpha=a_2=5$ and $\beta=b_1=1$, and so 
    \[O(G(S))-E(G(S)) = \dfrac{\alpha\beta-1}{2}=\dfrac{5\cdot1-1}{2}=2.\] 
    We note that $G(S)=\{1,3\}$, producing the same result for $O(G(S))-E(G(S))$ with direct calculation.
\end{example}

\section{Final thoughts}
\label{sec:final-thoughts}
\subsection{Alternate approaches with even generating elements}
Looking back, all of our results in this paper rely on an explicit description of an Ap\'ery set of a numerical semigroup. We have Proposition~\ref{prop:free-apery-description}, which gives the Ap\'ery set of a free numerical semigroup relative to the initial generating element $t_0$. If this generating element $t_0$ is odd, then we may use the function $f(n)=(-1)^n$ with Proposition~\ref{prop:shor-identity} to get a formula for $O(G(S))-E(G(S))$. And if this generating element $t_0$ is even, we need a different function $f(n)$. (Try $f(n)=(-1)^n$ to see why!) It turns out that $f(n)=n(-1)^n$ works and gets us our result.

In this subsection, we will briefly describe two alternate approaches to compute $O(G(S))-E(G(S))$ for $S$ a numerical semigroup generated by a compound sequence.

First, numerical semigroups generated by compound sequences have a few easily describable Ap\'ery sets. For instance, for $(A,B)$ a suitable pair and $C(A,B)=(n_0,\dots,n_p)$, there is an explicit description for the Ap\'ery set of $S=\langle C(A,B)\rangle$ relative to \emph{any} generating element $n_i$, as is given in \cite[Theorem 15]{KiersONeillPonomarenko16}:
\begin{equation}\label{eqn:apery-any-elt-compound}
    \Ap(S;n_i) = \left\{\sum\limits_{\substack{j=0 \\ j\ne i}}^p n_j x_j : 0\le x_j<b_{j+1} \text{ for } j<i, \text{ and } 0\le x_j<a_j \text{ for } j>i \right\}.
\end{equation}
Since some generating element $n_i$ is odd, we can use Theorem~\ref{prop:shor-identity} with $\Ap(S;n_i)$ and $f(n)=(-1)^n$ to get a formula for $O(G(S))-E(G(S))$.

We should also mention that we used this description of $\Ap(S;n_i)$ behind the scenes in Example~\ref{ex:first-ex}. In that case, the numerical semigroup is $S=\langle 6,15,20\rangle$, generated by the compound sequence $C(A,B)$ where $A=(2,3)$ and $B=(5,4)$. By the explicit description of $\Ap(S;n_i)$ in Eq.~\eqref{eqn:apery-any-elt-compound}, 
\[
    \Ap(S;6)
    =\{15x_1+20x_2 : 0\le x_1<2\text{ and } 0\le x_2<3\}
\]
and
\[
    \Ap(S;15)
    =\{6x_0+20x_2 : 0\le x_0<5\text{ and } 0\le x_2<3\}.
\]

Alternatively, if $C(A,B)=(n_0,\dots,n_p)$ is a compound sequence and $n_i$ is some odd generating element, then we can reverse the order in which the terms $n_0,\dots,n_i$ appear in the sequence (leaving the other terms where they are), to obtain a new sequence which generates the same numerical semigroup. This new sequence $T=(n_i,n_{i-1},\dots,n_0,n_{i+1},n_{i+2},\dots,n_p)$ is not necessarily compound, but it is telescopic! (See \cite[Proposition 3.10]{GassertShor17}.) Then, for $A=(a_1,\dots,a_p)$ and $B=(b_1,\dots,b_p)$, we have $c(T)=(b_i,\dots,b_1,a_{i+1},\dots,a_p)$. We now have a telescopic sequence with odd initial term, so we only need the result from Theorem~\ref{thm:main result} involving an odd initial term to compute $O(G(S))-E(G(S))$.

\subsection{Taking it further}
In this paper, we have considered the distribution of gaps among odd and even numbers, which means we have looked at gaps in the two congruence classes modulo 2. It could be interesting to see what sort of results one can obtain working modulo $m$ in general.

\section*{Acknowledgments}
We wish to thank Pieter Moree for the suggestion to look at \cite{Cho03062025}, which is how this paper came to be.

\printbibliography

\end{document}